\documentclass[12pt]{amsart}
\usepackage{amsmath,amsthm,epsfig,amssymb}
\usepackage{amsthm,amsfonts,amssymb,amscd}
\usepackage[all]{xy}
\usepackage{young}
\newtheorem{thm}[subsection]{Theorem}
\newtheorem{prop}[subsection]{Proposition}

\newtheorem{lem}[subsection]{Lemma}

\newtheorem{rem}[subsection]{Remark}
\theoremstyle{definition}
\newtheorem{Def}[subsection]{Definition}
\newtheorem{Not}[subsection]{Notation}

\newtheorem{proposition-definition}[subsection]{Proposition-Definition}

\begin{normalsize} \end{normalsize}

\newcommand{\CC}{{\mathbb C}}

\newcommand{\PP}{{\mathbb P}}

\newcommand{\WWW}{{\mathcal W}}
\newcommand{\AAA}{{\mathcal A}}
\newcommand{\OOO}{{\mathcal O}}

\numberwithin{equation}{section}
\renewcommand\square{\frame{\phantom{{\large x}}}}

\author{F. Laytimi}

\address{F. L.: Math\'ematiques - b\^{a}t. M2, Universit\'e Lille 1,
F-59655 Villeneuve d'Ascq Cedex, France}
\email{fatima.laytimi@math.univ-lille1.fr}

\author{W. Nahm}

\address{W. N.: Dublin Institute for Advanced Studies,
10 Burlington Road, Dublin 4, Ireland} \email{wnahm@stp.dias.ie}


\title{Generation by sections and  $k$-ampleness }

\begin{document}

\date{}

\begin{abstract}
In the article ``Submanifold of abelian varieties'', A.J. Sommese proved that
direct sum and tensor product  of two vector bundles $E$ and $F$  over a smooth
projective variety are $k$-ample if  $E$ and $F$ are $k$-ample and 
are  generated by global sections. Here we show that the latter condition is not needed.

\end{abstract}

\maketitle

\section{Introduction} \setcounter{page}{1}
Let $X$ be a smooth projective variety. 
The following definition was introduced by Sommese \cite{S1}.

\begin{Def}  A line bundle $L$ on $X$ is $k$-ample,
if for some $r>0, \, L^r $ is generated by sections and the fibers of the corresponding morphism
$$\phi:X \to \PP H^0(X,L^r) $$
have  dimensions less or equal to $k. $

A vector bundle $E$ is said to be $k$-ample if $\OOO_{\PP E}(1)$
is $k$-ample. 
\end{Def}

Note that $0$-ample is the same as ample.

\begin{Not}\label{Notgbs}
Since we often need to use the phrase ``generated by sections'', we will abbreviate it by ``gbs''.
\end{Not}

Sommese proved in ( \cite{S1} p. 235, Corollary (1.10)) that direct sums and tensor products of $k$-ample vector bundles 
are $k$-ample, if these vector bundles are generated by sections. The aim of the present article is to remove the latter 
restriction.

\section{Generation by sections}
We start  by giving several interpretations of the statement:
         $$ \OOO_{\PP E}(r)  \textrm {\ is gbs}.$$ 

1) \  Let $\pi:\PP E\rightarrow X$ be the projection. The fibre of $\pi$ over a point $x \in X$ is given by the 1-codimensional
subspaces $V$ of $E_x$, where $E_x$ is the fibre of $E$ over $x$. Points of $\PP E$ will be denoted by pairs $(x,V)$.
Let $V=v^\perp$ for $v \in E_x^*$. A section $s\in H^0(\PP E,\OOO_{\PP E}(r))$ generates the fibre of $\OOO_{\PP E}(r)$ 
over $(x,V)$ iff
$$\langle s(x,V), v^r\rangle\neq 0,$$
where $\langle\  , \ \rangle$ denotes the pairing of dual vector spaces.\\

2)\  Equivalently, let $\tilde s$ be the section in $H^0(X,S^rE)$ which corresponds to $s$ by the natural isomorphism
$$H^0(\PP E,\OOO_{\PP E}(r))\simeq H^0(X,S^rE).$$ 
Then $s$ generates the fibre of $\OOO_{\PP E}(r)$ over $(x,v^\perp)$ iff
$$\langle \tilde s(x), v^{\otimes r}\rangle\neq 0.$$

3)\ Given $x\in X$, the value $\tilde s(x)$ of each section $\tilde s \in H^0(X,S^rE)$ can be regarded as a homogeneous 
polynomial map of degree $r$ from $E_x^*$ to $\CC$. If the fibre of $\OOO_{\PP E}(r)$ over $v^\perp$ is not generated by 
any global section, then these polynomials have a common zero at $v$. If $\OOO_{\PP E}(r)$ is gbs, 
then 0 is the only common zero.\\

4)\ Since $v^{\otimes r}$ annihilates the kernel of the natural map $E^{\otimes r}\twoheadrightarrow S^rE$,
the fibre of $\OOO_{\PP E}(r)$ over $(x,v^\perp)$ is generated by a section iff
there is a section $\tilde s \in H^0(X,E^{\otimes r})$ such that
$$\langle \tilde s(x), v^{\otimes r}\rangle\neq 0.$$

If $\OOO_{\PP E}(r)$ is gbs, then $\OOO_{\PP E}(kr)$ is gbs for every positive integer $k$. All these facts will be
used in the following without reference.

\begin{lem} \label{sumgbs} Let $X$ be a complex space not necessarily compact and $E,F$ vector bundles on $X$. Suppose 
that $\OOO_{\PP E}(r)$ and $\OOO_{\PP F}(r)$ are gbs \rm (see Notation \ref{Notgbs}\rm ) for some $r>0$. Then
 $\OOO_{\PP(E\oplus F)}(r)$ is gbs. 
\end{lem}

Proof:
For each $x\in X$ and each non-zero $u\in (E_x^*\oplus F_x^*)$ we have to find a section
in $H^0(X,S^r(E\oplus F))$ which is not annihilated by $u^{\otimes r}$. Consider  $u=u_E+u_F$,
where $u_E,u_F$ are the components in $E_x^*, F_x^*$ resp. We may assume that $u_E\neq 0$. Since $\OOO_{\PP E}(r)$ 
is gbs, there is a section in $H^0(X,S^rE)$ which is not annihilated by $u_E^{\otimes r}$. The image of this section 
under the natural injection $H^0(X, S^rE)\hookrightarrow H^0(X, S^r(E\oplus F))$ is not annihilated by $u^{\otimes r}$.
$\hfill{\square}$

\begin{lem} \label{productgbs} Let $X$ be a complex space, not necessarily compact and $E,F$ vector bundles on 
$X$. Suppose that
 $\OOO_{\PP E}(r)$ and $\OOO_{\PP F}(r)$ are gbs (see Notation \ref{Notgbs}) for some $r>0$. Let $e={\rm rk} E$. Then
 there exists a positive integer $n$ depending only on $r$ and $e$ such that $\OOO_{\PP (E\otimes F)}(nr)$ is gbs.
\end{lem}

The proof of this lemma needs some preparations. We will use the terminology of decorated
oriented graphs (see Feynman diagrams in physics). Decoration means that to each vertex and to each arrow
another object is associated. The set of arrows of a graph $\gamma$ will be denoted by $A(\gamma)$. 
Any arrow $a$ is said to have a head $h(a)$ and a tail $t(a)$. A vertex which is the head of exactly $n$ arrows
is said to have indegree $n$, a vertex which is the tail of exactly $n$ arrows is said to have outdegree $n$.
The bidegree of a vertex is written as (indegree, outdegree).

\begin{Def} 
${}$

Let $V,W$ be vector spaces. Let $i=1,\ldots,n$, $\tilde v_i\in S^r V$, $\tilde w_i\in S^r W$, $u\in (V\otimes W)^*$. Then
$\Gamma_{V,W}(n,r,\tilde v_1,\ldots,\tilde v_n,\tilde w_1,\ldots,\tilde w_n,u)$ is the set of all decorated oriented graphs 
with $2n$ fixed vertices $\alpha_i, \beta_i,$ such that for each $i$ the vertex $\alpha_i$ has bidegree $(0,r)$
and is decorated by $\tilde v(\alpha_i)=\tilde v_i$, the vertex $\beta_i$ has bidegree $(r,0)$ and is decorated by 
$\tilde w(\beta_i)=\tilde w_i$, 
and such that all $rn$ arrows are decorated by $u$. We will omit the indices $V,W$ of $\Gamma_{V,W}$, when no
ambiguity arises.
\end{Def}

\begin{Def} 
 Let $\gamma\in \Gamma(n,1,v_1,\ldots,v_n,w_1,\ldots,w_n,u)$. The value $\vert\gamma\vert$ of $\gamma$ is defined by 
$$\vert\gamma\vert = \prod_{a\in A(\gamma)} u\big(v(h(a))\otimes w(t(a))\big).$$ 
\end{Def}

\begin{Def} 
 Let $\gamma\in \Gamma(n,r,\tilde v_1,\ldots,\tilde v_n,\tilde w_1,\ldots,\tilde w_n,u)$. Then the value $\vert\gamma\vert$ of 
$\gamma$ is defined by the following two conditions:
Firstly, $\vert\gamma\vert$ depends linearly on all $\tilde v_i$ and $\tilde w_i$. Secondly, if
$\tilde v_i = v_i^{\otimes r}$, $\tilde w_i = w_i^{\otimes r}$ for $i=1,\ldots,n$, $v_i\in V$,  $w_i\in W$ let the expanded 
graph
$$\gamma_{ex}\in \Gamma(rn,1,\underbrace{v_1,\ldots,v_1}_{r\ times},\ldots,\underbrace{v_n,\ldots,v_n}_{r\ times},
\underbrace{w_1,\ldots,w_1}_{r\ times},\ldots,\underbrace{w_n,\ldots,w_n}_{r\ times},u)$$
be a decorated oriented graph with vertices $\alpha_i^l, \beta_i^l,$ $i=1,\ldots,n$, $l=1,\ldots,r,$ such that
there is a bijection $\xi:A(\gamma_{ex})\rightarrow A(\gamma)$ with
$h(a)=\alpha_i$ if $h(\xi(a))=\alpha_i^l$, $t(a)=\beta_i$ if $t(\xi(a))=\beta_i^l$. Moreover, let
$\alpha_i^l$ be decorated by $v_i$, $\beta_i^l$ by $w_i$, for $i=1,\ldots,n$, $l=1,\ldots,r,$. Then 
$$\vert\gamma\vert = \vert\gamma_{ex}\vert.$$ 
\end{Def}

For later use, note that $\gamma_{ex}$ has a symmetry group of order $(r!)^{2n}s_\gamma^{-1}$, where $s_\gamma$ is 
the order of the group of vertex preserving symmetries of $\gamma$.\\

\begin{prop}\label{nu} There is a function $\nu:N \times N \rightarrow  N$ with the following property.
Consider finite dimensional vector spaces $V,W$ with $d={\rm dim} V$ and subspaces $A\subset S^rV$,
$B \subset S^rW$, such that the corresponding spaces of polynomial maps have 0 as only common zero. 
Let $n=\nu(r,d)$. Then for any non-trivial $u \in (V\otimes W)^*$ there is a positive integer $m$ with $m\leq n$,
elements $\tilde v_i\in A$, $\tilde w_i\in B$ for $i=1,\ldots,m$ and a decorated directed graph
$\gamma\in \Gamma_{V,W}(m,r,\tilde v_1,\ldots,\tilde v_m,\tilde w_1,\ldots,\tilde w_m,u)$
such that  $\vert\gamma\vert\neq 0$.
\end{prop}

Proof: 
We will construct a function  $\nu$ which is far from optimal, but  sufficient for our purpose.

We first show that it suffices to prove the proposition for the case that the map $\hat u\in Hom(V,W^*)$ corresponding to 
$u$ is bijective. Otherwise, let $V'=V/ker(\hat u)$,  $W'=W/ker(\hat u^*)$, let $\pi_V,\pi_W$ be the corresponding projections 
$\pi_V:S^rV\twoheadrightarrow S^rV'$, $\pi_W:S^rW\twoheadrightarrow S^rW'$ and let $u'$ be the element of $(V'\otimes W')^*$
induced by $u$. If $$\gamma\in \Gamma(m,r,V_1,\ldots,V_m,W_1,\ldots,W_m,u), \rm {and} $$ 
$$\gamma'\in \Gamma(m,r,\pi_V(V_1),\ldots,\pi_V(V_m),\pi_W(W_1),\ldots,\pi_W(W_m),u')$$ 
have the same underlying undecorated graph, then $\vert\gamma'\vert=\vert\gamma\vert$.

If $u$ is bijective, we can identify $W$ with $V^*$ and write
$$u(v\otimes w)=\langle v,w\rangle$$  
for $v\in V$, $w\in W$. Let $SV$ be the symmetric algebra over $V$.
There are natural multiplication maps  $$m^r:S^rV\rightarrow End(SV),$$ 
and contraction maps  $$i^r:S^rV \rightarrow End(SV)$$
which restrict to
$$m^r_N:S^rV\rightarrow Hom(S^{N-r}V,S^NV),$$ 
$$i^r_N:S^rW \rightarrow Hom(S^NV,S^{N-r}V),$$ 
with $S^NV=0$ for negative $N.$ For $r=1$ the maps can be characterized by multilinearity and the properties
$m^1(v')v^{\otimes N}=\pi_S (v'\otimes v^{\otimes N})$, $i^1(w)v^{\otimes N}=N\langle w,v\rangle v^{\otimes (N-1)}$
for $v',v\in V$, $w\in W$, where $\pi_S$ is the projection of the tensor algebra of $V$ to the symmetric algebra $SV$.
For general $r$ they are characterized by multilinearity and the properties $m^r(v^{\otimes r})=m^1(v)^r$, 
$i^r(w^{\otimes r})=i^1(w)^r$.

For $r=1$ let $P$ be a product with the factors $m^1(v_i)$, $i^1(w_i)$, taken in any order, where $v_i\in V$, $w_i\in W$,
$i=1,\ldots,m$. For any non-negative integer $N$, the product $P$ restricts to a map $P_N\in End(S^NV)$. For the trace of $P_N$ 
one finds
$$tr (P_N) = \sum_{\gamma\in \Gamma_m} c_\gamma\vert\gamma\vert,$$
where 
$$\Gamma_m = \Gamma_{V,W}(m,1,v_1,\ldots,v_m,w_1,\ldots,w_m,\langle,\rangle)$$
and 
$$c_\gamma={d+\rho+N-1\choose d+m-1},$$
where $\rho$ is the cardinality of the set of arrows $a$ of $\gamma$ such that the factor $m^1(v(h(a)))$ lies
to the right of the factor  $i^1(w(t(a)))$ in $P$. For $\rho>0$ this follows by induction on $\rho$
from
$i^1(w)m^1(v)=m^1(v)i^1(w)+\langle w,v\rangle$ and for $\rho=0$ from cyclic invariance of the trace and induction on $N$.
For $n=0$ and $N=0$ the statement is obvious. 

For arbitrary $r$ the calculation can be reduced to the case $r=1$ as in the definition of $\vert\gamma\vert$. 
Let  $P$ be a product with the factors $m^r(\tilde v_i)$, $i^r(\tilde w_i)$, taken in any order, 
where $\tilde v_i\in S^rV$, $\tilde w_i\in S^rW$, $i=1,\ldots,m$. For any $N$, $P$ restricts to a map $P_N\in End(S^NV)$. 
To calculate the trace of $P_N$ it is sufficient to consider the case $\tilde v_i = v_i^{\otimes r}$, 
$\tilde w_i = w_i^{\otimes r}$ for $i=1,\ldots,n$, $v_i\in V$,  $w_i\in W$. One finds
$$tr( P_N) = \sum_{\gamma\in \Gamma^r_m} c_\gamma\vert\gamma\vert,$$
where 
$$\Gamma^r_m = \Gamma_{V,W}(m,r,\tilde v_1,\ldots,\tilde v_m,\tilde w_1,\ldots,\tilde w_m,\langle,\rangle)$$
and 
$$c_\gamma=(r!)^{2m}s_\gamma^{-1}{d+r\rho+N-1\choose d+rm-1}.$$
Here $\rho$ is the cardinality of the set of arrows $a$ of $\gamma$ such that the factor $m^r(\tilde v(h(a)))$ lies
to the right of the factor  $i^r(\tilde w(t(a)))$ in $P$, and  $s_\gamma$ is the order of the group of vertex preserving
symmetries of $\gamma$.

Let $\AAA^r_N$ be the subalgebra of ${\rm End}(S^NV)$ generated by the elements
$\{m^r_N(\tilde v) i^r_N(\tilde w)\vert \tilde v\in A, \tilde w\in B\} $. This algebra is spanned by products
$P_N=m^r_N(\tilde v_1) i^r_N(\tilde w_1)\cdots m^r_N(\tilde v_m) i^r_N(\tilde w_m)$, where $\tilde v_i\in A$,
$\tilde w_i\in B$ for $i=1,\ldots,m$.
We have seen that the trace of $P_N$ is given by a linear combination of numbers
$\vert\gamma\vert$, where $\gamma\in \Gamma(m,r,\tilde v_1,\ldots,\tilde v_m,\tilde w_1,\ldots,\tilde w_m,\langle,\rangle)$.
If those traces all vanish, then $ \AAA^N$  is a nilalgebra. In particular, there is a non-trivial subspace $V_0\subset S^NV$
such that 
$$m^r_N(\tilde v) i^r_N(\tilde w)V_0=0$$ 
for all $\tilde v\in A$, $\tilde w\in B$.
Since the kernel of $m^r_N(\tilde v)$ vanishes for $\tilde v\neq 0$, this means that 
$i^r_N(\tilde w)V_0=0$ for all $\tilde w\in B$. 
Since the homogeneous polynomial maps in $B$ have 0 as only common zero, one can apply
a theorem of Macaulay (\cite{M} Sec. 6, or \cite{S2} p. 85 theorem 4.48). According to this theorem the ideal 
in $SW$ generated by $B$ contains $S^NW$ if $N\geq rd$, which implies $i^N_N(\tilde w)V_0=0$
 for all $\tilde w\in S^NW$ and
yields a contradiction. Thus $\AAA^N$ is not nil for $N\geq rd$.

Let $\WWW^N_m$ be the subspace of $\AAA^N$ generated by products of length $\leq m$ of elements 
of the form $m^r_N(\tilde v) i^r_N(\tilde w)$.
For each $m$ one has either\\  ${\rm dim} (\WWW^N_{m+1})>{\rm dim} (\WWW^N_m$) or $\WWW^N_m = \AAA^N$.\\  Let 
$D(N)=({\rm dim} \ S^N V)^2$. Since $D(N)\geq {\rm dim} \AAA^N$, one has $\WWW^N_{D(N)} = \AAA^N$. In particular,
$\WWW^{rd}_{D(rd)}$ contains elements of non-vanishing trace. Thus it suffices to take for $\nu(r,d)$ the least 
common multiple of all integers less or equal to $D(rd)$.

 $\hfill{\square}$\\

{\bf Proof of Lemma}  \ref{productgbs}:\\

For each $x\in X$ and each non-zero $u\in (E_x\times F_x)^*$ we have to find a section
in $H^0(X,(E\otimes F)^{\otimes nr})$ which is not annihilated by $u^{\otimes nr}$.
Let 
$$I_E: H^0(X,S^rE)^{\otimes n}\hookrightarrow  H^0(X,E^{\otimes rn}),$$
$$I_F: H^0(X,S^rF)^{\otimes n}\hookrightarrow  H^0(X,F^{\otimes rn}),$$
$$J: H^0(X,E^{\otimes rn})\otimes H^0(X,F^{\otimes rn})\hookrightarrow H^0(X,(E\otimes F)^{\otimes rn})$$
be the canonical injections.
The permutation group $S(rn)$ acts on $E^{\otimes rn}$ in the standard way, which yields a map
$$\Sigma:S(rn)\times H^0(X,E^{\otimes rn}) \rightarrow  H^0(X,E^{\otimes rn}).$$
Let
$$\Phi:S(rn)\times H^0(X,S^rE)^{\otimes n}\otimes H^0(X,S^rF)^{\otimes n}\rightarrow H^0(X,(E\otimes F)^{\otimes nr})$$
be defined by $\Phi=J\circ(\Sigma\otimes Id)\circ(Id\times I_E\otimes I_F)$.
Though $\Phi$ factors through the coset map $S(rn)\rightarrow S(rn)/S(r)^{\times n}$, it yields sufficiently
many sections for our purpose.
For $x\in X$, $\sigma\in S(rn)$, $\tilde s_i\in H^0(X,S^rE)$, $\tilde t_i\in H^0(X,S^rF)$, there is an element
$$\gamma\in \Gamma_{E_x,F_x}(n,r,\tilde s_1(x),\ldots,\tilde s_n(x),\tilde t_1(x),\ldots,\tilde t_n(x),u),$$ 
such that

\begin{equation} \label{phi}
\langle\Phi(\sigma,\tilde s_1\otimes \cdots \otimes\tilde s_n\otimes \tilde t_1\otimes \cdots \otimes\tilde t_n)_x,
 u^{\otimes nr}\rangle=\vert\gamma\vert.
\end{equation}

Conversely, for each $\gamma\in \Gamma(n,r,\tilde s_1(x),\ldots,\tilde s_N(x),\tilde t_1(x),\ldots,\tilde t_n(x),u)$
one can find a permutation $\sigma\in S(rn)$ satisfying  the equation\eqref{phi}.

By Proposition \ref{nu} there exist sections $\tilde s_i$, $\tilde t_i$ and 
$$\gamma\in \Gamma(n,r,\tilde s_1(x),\ldots,\tilde s_N(x),\tilde t_1(x),\ldots,\tilde t_n(x),u)$$ 
such that $\vert\gamma\vert\neq 0$. The image of the corresponding section 
$\Phi(\sigma,\tilde s_1\otimes \cdots \otimes\tilde s_n\otimes \tilde t_1\otimes \cdots \otimes\tilde t_n)$
in $H^0(X,S^{\otimes nr}(E\otimes F))$ yields a section in  $H^0(\PP (E\otimes F),\OOO_{\PP (E\otimes F)}(nr))$
which generates the fibre over $(x,u^\perp)$.
 $\hfill{\square}$

\section {Direct sums and tensor products of $k$-ample vector bundles}

\begin{thm}\label{sum}

If $E,F$ are  $k$-ample, then $E\oplus F$ is  $k$-ample.
\end{thm}

Proof: We use criterion (1.7.3) of proposition (1.7) of Sommese \cite{S1}.
First note that $\OOO_{\PP E }(r_1)$ gbs (see Notation \ref{Notgbs}) and $\OOO_{\PP F}(r_2)$ gbs imply
that $\OOO_{\PP E }(r)$ and $\OOO_{\PP F }(r)$ are gbs whenever $r$ is a common multiple of $r_1,r_2$. 
By Lemma  \ref{sumgbs} this implies that $\OOO_{\PP (E\oplus F)}(r)$ is gbs. 
Now assume that there is a holomorphic finite to one map $\phi:Z\rightarrow X$
of a compact analytic space $Z$ to $X$, such that ${\rm dim} Z = k+1$ and that there is a surjective map 
$\phi^*(E\oplus F)\twoheadrightarrow Q$
with a trivial bundle $Q$. In particular there is a non-trivial section of
$\phi^*(E^*\oplus F^*)$. We may assume that the component $s:Z \rightarrow \phi^* E^*$
of this section is non-trivial.
Since $\OOO_{\PP E}(r)$ is gbs, for any $z\in Z$ such that $s(z)\neq 0$ there is a section
$\sigma: X\rightarrow S^rE$ such that $\langle \phi^*\sigma(z), s(z)^{\otimes r}\rangle\neq 0$.
In particular, $\langle \phi^*\sigma, s^r\rangle$ yields a non-trivial section
of the trivial line bundle over $Z$. Such a section must be constant, which implies
that $s$ cannot vanish anywhere. Consequently it yields a trivial quotient bundle of $\phi^*E$,
contrary to the assumption that $E$ is $k$-ample. $\hfill{\square}$

\begin{thm}\label{product}

If $E,F$ are $k$-ample, then $E\otimes F$  is $k$-ample.
\end{thm}

Proof: We use criterion (1.7.4), proposition (1.7)  of Sommese \cite{S1}. We first show that $S^2E$ is $k$ ample.
By Lemma 1 the bundle $E^{\oplus e}$ is $k$-ample, where $e$ is the rank of $E$. 
Since each irreducible summand of $S^n(S^2E)$ is isomorphic to an irreducible summand of $S^{2n}(E^{\oplus e}),$  
we have $H^j(X,F\otimes S^n(S^2E))=0$ for $j>k$, for any coherent sheaf $F$ on $X$ and $n>>0$. 
Since $\OOO_{\PP (S^2E)}(1)$ is a restriction of $\OOO_{\PP (E^{\otimes 2})}(1)$, by Lemma \ref{productgbs}
$\OOO_{\PP (S^2E)}(N)$ is gbs for some N. Thus $S^2E$ is $k$ ample. 
Since $E\otimes F$ is a quotient of the $k$-ample bundle $S^2(E\oplus F)$, the theorem follows. 
$\hfill{\square}$

\begin{rem}
 Theorem \ref{product} implies that there is no need for the gbs conditions in Lemma (1.11.4) and Proposition (1.13) of
\cite{S1}. One obtains the following.
\end{rem}

\begin{lem}
 Let $E$ be a vector bundle on a compact analytic space $X$, $Gr(s,E)$ be the bundle of $s$-codimensional subspaces
of the fibres of $E$ and $\xi(s,E)$ the tautological line bundle of $Gr(s,E)$. Then $\xi(s,E)$ is $k$-ample if $E$ is
$k$-ample.
\end{lem}

\begin{prop} For $i=1,\ldots,n$ let $E_i$ be a vector bundle of rank $r_i$ on a projective manifold $X$, If for all $i$ 
$E_i$ is $k_i$-ample. Then
$$H^p(X,\wedge^qT^*_X\otimes \wedge^{s_1}E_1\otimes \cdots \otimes \wedge^{s_n}E_n)=0$$
if $p+q>{\rm dim} X + \sum_i s_i(r_i-s_i) + {\rm min}_i\{k_i\}.$
\end{prop}


\begin{thebibliography}{000}

\bibitem{M} P. Macaulay, {\it The algebraic geometry of modular systems, }
Cambridge Univ. Press, 1916.

\bibitem{S2} B. Shiffmann, A.J. Sommese, {\it Vanishing theorems on complex manifolds, }
Progress in mathematics Vol. {\bf 56}(1985).

\bibitem{S1} A.J. Sommese, {\it Submanifold of Abelian
Varieties,}, Math.Ann.{\bf 233}(1978), 229-256.


\end{thebibliography}
\end{document}